\documentclass[11pt]{article}

\usepackage{amsmath,amssymb,amsthm, color}
\newtheorem{theorem}{Theorem}[section]
\newtheorem{lemma}[theorem]{Lemma}
\newtheorem{proposition}[theorem]{Proposition}
\newtheorem{definition}[theorem]{Definition}

\newtheorem{remark}[theorem]{Remark}
\numberwithin{equation}{section}

\def\({\bigl(}                \def\){\bigr)}

                        \def\^{\tilde}

                      \def\R{I\!\!R}

\def\1{1\!\!1}

\date{\today}

\begin{document}
\title{The extended Hausdorff dimension spectrum of a conformal iterated function system is maximal}

\author{
Andrei E. Ghenciu
\thanks{The first author would like to thank
the University of Denver for its hospitality during
his stays while part of this paper was written.},\\
Department of Mathematics, Statistics and Computer Science,\\
University of Wisconsin Stout, Menomonie, WI, 54751, USA.\\
email: ghenciue@uwstout.edu
\\ \\and\\ \\
Ronnie Pavlov
\thanks{The second author gratefully acknowledges the support of a Simons collaboration grant.},\\
Department of Mathematics,\\
University of Denver, \\
2390 S. York St., Denver, CO, 80210, USA.\\
email: rpavlov@du.edu
}
\date{\today}

\maketitle

\begin{abstract}
For any conformal iterated function system (CIFS) consisting of finitely or countably many maps, and any closed shift-invariant set of right-infinite sequences of such maps, one can associate a limit set, which we call a shift-generated conformal iterated construction. We define the extended Hausdorff dimension spectrum of a CIFS to be the set of Hausdorff dimensions of all such limit sets. We prove that 
for any CIFS with finitely or countably many maps, the extended Hausdorff dimension spectrum is maximal, i.e. all nonnegative dimensions less than or equal to the dimension of the limit set of the CIFS are realized. We also prove a version of this result even for so-called conformal graph directed Markov systems, obtained via nearest-neighbor restrictions on the CIFS.

The main step of the proof is to show that for the family $(X_\beta)$ of so-called $\beta$-shifts, the Hausdorff dimension of the limit set associated to $X_\beta$ varies continuously as a function of $\beta$.
\end{abstract}
\section{Introduction}
The study of iterated function systems (IFSes) began in the early 1980s. Many leading researchers made significant contributions to the theory of IFSes and various extensions, including Barnsley, Falconer, Graf, Hutchinson, Mauldin, and Urbanski. 

The initial focus was on IFSes consisting of finitely many Euclidean similarities. Later on, the theory was extended to allow infinitely many maps, often called infinite IFSes. Mauldin and Urbanski were among the first to pioneer two extensions of IFS theory, first the infinite conformal iterated function systems (CIFSes) (see \cite{MU}), and their generalizations, conformal graph directed
Markov systems (CGDMSes) (see \cite{MU4}). We postpone formal definitions to Section~\ref{prelim}, but any CIFS/CGDMS is defined by a countable family $(\phi_a)_{a \in A}$ of contractions on a compact metric space $Y$. A CIFS allows for arbitrary composition of the maps, whereas a CGDMS has nearest-neighbor restrictions on how the maps can be applied; for instance perhaps $\phi_2$ cannot be applied after $\phi_5$.

A well-studied question in this area (first asked by Mauldin and Urbanski in \cite{mutr}) is: given a CIFS with attractor $J$, what are the possible Hausdorff dimensions of limit sets of sub-CIFSes? They defined the \textbf{Hausdorff dimension spectrum} of a CIFS 
with (infinite) alphabet $A$ as the set $\{ HD(J(B)) : B \subset A\}$, where $J(B)$ is the limit set of the restricted CIFS with maps 
$(\phi_b)_{b \in B}$, and $HD(S)$ is the Hausdorff dimension of a set $S$. Clearly this set is a subset of $[0, HD(J)]$, but how large is it?

In the case of the standard continued fractions CIFS over the unit interval (indexed by $A = \mathbb{N}$), an extremely similar question was asked independently by Hensley (\cite{He}): must the set $\{ HD(J(B)) : B \subset \mathbb{N}, |B| < \infty\}$ be dense in 
$[0,1]$? This was dubbed the Texan conjecture by Jenkinson, and was resolved in the affirmative in 2006 by Kessebohmer and Zhu in \cite{KZ}.
Research surrounding the Texan conjecture gave birth to techniques for studying topological and metric features of the Hausdorff dimension spectrum for a general CIFS.

Mauldin and Urbanski showed that the closed interval $[0, \theta]$ of any CIFS is always contained in the Hausdorff dimension 
spectrum, where $\theta$ is the so-called finiteness parameter of the CIFS.

The Hausdorff dimension spectrum has been studied almost exclusively in the case of an infinite alphabet, because when $A$ is finite, there are only finitely many subsets $B \subset A$, and so the traditionally defined Hausdorff dimension spectrum is uninteresting.

However, in analogy to the way that a CDGMS is defined by restricting the possible `itineraries' of applications of the maps $\phi_a$ to a topological Markov shift, one can consider much more general subsystems of a CIFS. Rather than removing maps from the collection entirely, one can restrict the sequences of legal applications to an arbitrary shift space $X$ on alphabet $A$; we call this a
\textbf{shift-generated conformal iterated construction}. Even when $A$ is finite, there are uncountably many such subsystems, and each yields a well-defined limit set $J(X)$, and so the collection of their Hausdorff dimensions (which we call the \textbf{extended Hausdorff dimension spectrum}) could potentially be quite large.

In fact this is our main result: for any CIFS or CDGMS over a finite alphabet $A$, the extended Hausdorff dimension spectrum is
$[0, HD(J)]$. The main idea of our proof is to use the so-called $\beta$-shifts, a one-parameter collection $(X_{\beta})$ of subshifts over any finite alphabet $\{0, \ldots, k\}$, and to prove that $HD(J(X_{\beta}))$ is a continuous function of $\beta$. The utility of 
$\beta$-shifts here is not a coincidence; they are one of the simplest families of subshifts known to achieve any possible (topological) entropy (see \cite {Wa}, for more details), which is related to the Hausdorff dimension of limit sets via the Bowen formula (presented here as Proposition~\ref{Bowenf}).


\section{Preliminaries}\label{prelim}

\subsection{Symbolic dynamics}

We begin with some definitions, results and concepts from symbolic dynamics; for further background and proofs, we refer the reader to \cite{LM}. 

Define $A$ to be any finite or countable nonempty set, called the \textbf{alphabet}. The elements of $A$ will be called \textbf{letters}.
\begin{definition}
A \textbf{word} over $A$ is any finite sequence of letters from $A$. More formally, a word is any $\omega \in A^n$, for some $n \geq 1$; this $n$ is referred to as the \textbf{length} of $\omega$.
\end{definition}
\begin{definition}
The \textbf{(one-sided) full $A$-shift} is the set of all infinite sequences of letters from $A$:
$$ A^{\mathbb{N}}= \{ x = (x_1, x_2, ..., x_i,...), x_i \in A \}. $$
\end{definition}
\begin{definition}
The \textbf{(left) shift map} $\sigma$ from $A^{\mathbb{N}}$ is defined by $(\sigma x)_n = x_{n+1}$ for all $n \geq 0$, i.e.
\[
(x_1, x_2, \ldots) \mapsto (x_2, x_3, \ldots).
\]
\end{definition}
\begin{definition}
A \textbf{shift space} or \textbf{subshift} over the alphabet $A$ is a set $X \subset A^{\mathbb{N}}$ which is closed (in the product discrete topology) and invariant under $\sigma$, i.e. $x \in X \Longrightarrow \sigma x \in X$.
\end{definition}
\begin{definition}
A \textbf{nearest-neighbor shift of finite type} or \textbf{topological Markov chain} is a subshift over the alphabet $A$ defined by a set $S \subset A^2$ of adjacencies in the following way: 
\[
X_S = \{(x_1, x_2, \ldots) \in A^{\mathbb{N}} \ : \ \forall n, (x_n, x_{n+1}) \in S\}.
\]
A topological Markov chain over a finite alphabet is called \textbf{irreducible} if for all $i, j \in A$, there exists $w \in \mathcal{L}(X_S)$ beginning with $i$ and ending with $j$. 
A topological Markov chain over a countable alphabet is called \textbf{finitely irreducible} if there exists $N$ so that for all $i, j \in A$, there exists $w \in \mathcal{L}(X_S)$ with length at most $N$ beginning with $i$ and ending with $j$.
\end{definition}
\begin{definition}
For a subshift $X$ over $A$, the \textbf{language} of $X$ is the set $\mathcal{L}(X)$ of all words occurring within points of $x$, i.e. 
\[
\mathcal{L}(X) = \{x_k \ldots x_\ell \ : \ x \in X, k \leq \ell \in \mathbb{N}\}.
\]
For any $n$, the \textbf{$n$-language} of $X$ is $\mathcal{L}_n(X) := \mathcal{L}(X) \cap A^n$, the set of $n$-letter words in $\mathcal{L}(X)$. 
\end{definition}
For every $\omega\in{\mathcal L}(X)\cup X$ and every $n \in \mathbb{N}$ less than $|\omega|$ (which is taken to be $\infty$ if $\omega \in X$), we denote by $\omega|_n$ the word $\omega_1\omega_2\ldots\omega_n$. 
\subsection{$\beta$-shifts}

For $\beta > 1$, the $\beta$-shift $X_\beta$, which is a subshift with alphabet $\{0, 1, \ldots, \lceil \beta \rceil - 1\}$, was originally defined by Parry as a symbolic coding of the linear map $x \mapsto \beta x \mod 1$. The definition we give here uses greedy expansions with base $\beta$. For any $t \in [0,1)$, define the sequence $a(t)$ to be the greedy expansion of $t$ in powers of 
$\beta^{-1}$: $a_1$ is the maximal integer for which $a_1 \beta^{-1} < t$, and then for all $k > 1$, given $a_1, \ldots, a_k$, we define $a_{k+1}$ as the maximal integer for which $a_{k+1} \beta^{-(k+1)} < t - \sum_{i = 1}^k a_i \beta^{-i}$.

Then, the subshift $X_\beta$ can be defined as the closure of the set of all such sequences $a(t)$. It is easily checked that a sequence $x$ is in $X_\beta$ iff for all $0 \leq i < j$, $\sum_{s = 1}^{j-i} \frac{x_{s+i}}{\beta^s} < 1$.




\subsection{Shift-generated conformal iterated constructions over a finite or countable alphabet}


A \textbf{shift-generated iterated construction} is defined by 
a subshift $X$ over an alphabet $A$ with at least two letters, a non-empty compact metric space $Y$, 
a finite collection of non-empty compact subsets $\{Y_e\}_{e \in A}$ of $Y$, and 
a set of generators $\Phi=\{\phi_e:Y_e\to Y\}_{e\in A}$, where the $\phi_e$'s are one-to-one 
contractions which satisfy $\phi_f(Y_f)\subseteq Y_e$ 
whenever $ef\in{\mathcal L}_2(X)$. 

Let $0<s<1$ be such that all these generators 
have a contraction ratio that does not exceed $s$. For every 
$\omega\in{\mathcal L}(X)$, set $Y_\omega=Y_{\omega_{|\omega|}}$ and
\[
\phi_\omega:Y_\omega\to Y,\ \ \ \ 
\phi_\omega:=\phi_{\omega_1}\circ\phi_{\omega_2}\circ\ldots\circ\phi_{\omega_{|\omega|}}.
\]

Given $\omega\in X$, the compact sets $\phi_{\omega|_n}(Y_{\omega|_n})$, $n\geq1$, are decreasing
and their diameters converge to zero. More precisely,
\[
\mbox{diam}\bigl(\phi_{\omega|_n}(Y_{\omega|_n})\bigr)\leq s^n\mbox{diam}(Y).
\]

This implies that the set
\[
\bigcap_{n\geq1}\phi_{\omega|_n}(Y_{\omega|_n})
\]
is a singleton. We define the \textbf{coding map} $\pi:X\to Y$ by
\[
\{\pi(\omega)\}=\bigcap_{n\geq1}\phi_{\omega|_n}(Y_{\omega|_n})
\]
and we define the \text{limit set of the shift-generated construction} to be
\[
J(Y)=\pi(X).
\]

We call a shift-generated iterated construction \textbf{conformal} 
if the following conditions are satisfied:
\\ \\
(i)\ \ There exists $d$ so that for all $e \in A$, $Y_e$ is a connected compact subset of $\mathbb{R}^d$ and 
$Y_e=\overline{\mbox{Int}_{\mathbb{R}^d}(Y_e)}$.
\\ \\
(ii)\ \ (Open Set Condition (OSC)) For every $e,f\in A$, $e\ne f$,
\[
\phi_e(\mbox{Int}(Y_e))\cap\phi_f(\mbox{Int}(Y_f))=\emptyset.
\]
(iii)\ \ For every $f\in A$, there exists a connected open set $W_f$ with $Y_f\subseteq W_f\subseteq\R^d$
so that the map $\phi_f$ extends to a $C^1$ conformal
diffeomorphism of $W_f$ into $\displaystyle\bigcap_{e\in A:ef\in{\cal L}_2(X)}W_e$.
\\ \\
(iv)\ \ There are two constants $L\geq1$ and $\alpha>0$ so that
\[
\bigl||\phi_e'(x)|-|\phi_e'(y)|\bigr|\leq L\|(\phi_e')^{-1}\|^{-1}\cdot|x-y|^{\alpha}
\]
for every $e\in A$ and for every pair of points $x,y\in Y_e$, where
$|\phi_e'(x)|$ represents the norm of the derivative. 
\\ \\
(v)\ \ (Cone Property) There exists $\gamma,l>0$, $\gamma<\frac{\pi}{2}$ such that for every
$x\in Y$ there exists an open cone $\mbox{Con}(x,\gamma,l)\subseteq\mbox{Int}(Y)$ with vertex $x$, 
central angle $\gamma$, and altitude $l$.

\begin{remark}
If $d \geq 2$ and a construction satisfies conditions~(i) and~(iii), then it also 
satisfies condition~(iv) with $\alpha=1$ according to Proposition~4.2.1 in~\cite{MU4}.
\end{remark}

As a consequence of~(iv) we get the following:\\
\\
(iv')\ \ (Bounded Distortion Property (BDP)) There exists $K\geq1$ such that for all 
$\omega\in{\mathcal L}(X)$ and for all $x,y\in W_\omega$,
\[
|\phi_\omega'(y)|\leq K|\phi_\omega'(x)|.
\]

We note that in using our terminology, a \textbf{conformal iterated function system (CIFS)} is a shift-generated conformal iterated construction induced by a full shift (and $Y_e = Y$, for every $e \in A$), and a \textbf{conformal graph directed Markov system (CDGMS)} is a shift-generated conformal iterated construction induced by a topological Markov chain.

Next, we collect necessary versions of Bowen's formula for shift-generated conformal iterated constructions over a finite alphabet.
We start by defining one of the main tool in studying our constructions, a topological pressure - like function. For simplicity, we will call it the topological pressure function.
Given $t\geq0$ and $n\geq1$, we denote the $n$th-level partition function $Z(Y,n,t)$ by 
\[
Z(Y,n,t)=\sum_{\omega\in{\mathcal L}_n(X)}\|(\phi_\omega)'\|^{t}.
\]
For every $t\geq0$, the sequence $(Z(Y,n,t))_{n\geq1}$ is submultiplicative and thus we can define the 
topological pressure function $P(Y,t)$ of the construction by:
\[
P(Y,t)=\lim_{n\to\infty}\frac{1}{n}\log Z(Y,n,t)=\inf_{n\geq1}\frac{1}{n}\log Z(Y,n,t).
\]
The topological pressure function $P(Y,t):[0,\infty)\to\R$ (as a function of $t$) is strictly decreasing to negative infinity, convex and hence 
continuous. Indeed, the strictly decreasing behavior of the pressure can be more precisely described as follows.
Let $0\leq t_1<t_2$. Then $Z(Y,n, t_2)\leq s^{n(t_2-t_1)}Z(Y,n, t_1)$ for all $n\geq1$. Therefore
$P(Y,t_2)\leq(t_2-t_1)\log s+P(Y,t_1)$. The convexity of the pressure follows from the convexity 
of its partition functions $Z(Y,n)$. The continuity is a direct consequence of the convexity. \\
\\
One of the most important results for shift-generated conformal iterated constructions over a finite alphabet connects the Hausdorff dimension of the limit set with the zero of the topological pressure function.
\begin{proposition}[Bowen formula]\label{Bowenf}
Let $\Phi$ be a shift-generated conformal iterated construction over a finite alphabet, with subshift $X$ and underlying space $Y$, and let $h$ be 
the zero of its topological pressure function. Then $\mbox{HD}(J(Y))=h$ and 
$\mathcal{H}^h (J(Y))>0$, where $\mathcal {H}^h$ 
stands for the Hausdorff measure in the dimension $h$.
\end{proposition}
For a detailed proof of this result, please see \cite{GR}.
For general definitions and results related to the Hausdorff dimension and Hausdorff measure, see \cite{Fal}.

\section{Main Results}
The following lemma is the main technical tool in the proofs of our main results.
\begin{lemma}\label{beta}
 For any $\beta > 1$, $k \in \mathbb{N}$, $\delta < (1 + \beta^{2k})^{1/2k} - \beta$, and $y \in X_{\beta + \delta}$, there exists a set $S \subset \mathbb{N}$ satisfying the following:\\

\noindent
(1) $S$ has gaps greater than $k$, i.e. $s \neq t \in S \Longrightarrow |s - t| > k$

\noindent 
(2) $y(s) > 0$ for all $s \in S$

\noindent
(3) The sequence $x$ defined by replacing all letters of $y$ at positions in $S$ by $0$s is in $X_\beta$.

\end{lemma}

\begin{proof}

Choose any such $k, \beta, \delta, y$ and note that $(\beta + \delta)^{2k} < 1 + \beta^{2k}$. Recall our definition of $X_{\beta}$: $x \in X_{\beta}$ iff for every $0 \leq i < j$,
\begin{equation}\label{beta2}
\sum_{s = 1}^{j-i} \frac{x_{s+i}}{\beta^s} < 1.
\end{equation}


Now, assume that $y \in X_{\beta + \delta}$. We may partition $y$ into maximal runs of $0$s of length at least $k$ and remaining intervals, i.e. either $y = w_1 r_1 w_2 r_2 w_3 r_3 \ldots$ or $y = w_1 r_1 w_2 r_2 \ldots w_n r_n$ ($r_n$ infinite) or 
$y = w_1 r_1 \ldots r_{n-1} w_{n}$ ($w_{n}$ infinite) where each $r_i = 0^{n_i}$ for $n_i \geq k$ and each finite $w_i$ ends with a non-$0$ letter and contains no $0^k$. For all $i$ with $w_i$ finite (say $y([a_i, b_i]) = w_i$), define elements of $[a_i, b_i]$ as follows:
$s_i^{(1)} = b_i$ (note that $y\left(s_i^{(1)}\right) > 0$) and, for all $m$, if $s_i^{(m)}$ is defined and at least $a_i + k$, take
$s_i^{(m+1)} = \sup\{j \ : \ y(j) > 0, a_i \leq j < s_i^{(m)} - k\}$ (the set inside the supremum is nonempty since $w_i$ does not contain $k$ consecutive $0$s). If $s_i^{(m)} < a_i + k$, then stop and define $S_i = \{s_i^{(1)}, \ldots, s_i^{(m)}\}$. If $y$ terminates with an infinite $w_n$ (say $y([a_n, \infty)) = w_n$, then define elements of $[a_n, \infty)$ as follows: $s_n^{(1)} = a_n$ (note that $y\left(s_n^{(1)}\right) > 0$) and, for all $m$, take
$s_n^{(m+1)} = \inf\{j \ : \ y(j) > 0, j > s_n^{(m)} + k\}$ (the set inside the infimum is nonempty for the same reason as above). Note for future reference that, since every $w_i$ contains no $0^k$, 
$s_i^{(m+1)} \leq s_i^{(m)} + 2k$ for all $m,i$.

Define $S$ to be the set of all $s_i^{(j)}$ generated above. It is immediate from definition that $S$ satisfies (1) and (2). We need to verify that $x$ defined by replacing all letters in $y$ at locations in $S$ by $0$ is in $X_\beta$.

It suffices to check (\ref{beta2}) for $0 \leq i < j$ satisfying $x_{i+1} > 0$. 
This is because if 
$\sum_{s = 1}^{j-i} \frac{x_{s+i}}{\beta^s} = K \geq 1$ and $x_{i+1} = 0$, then we can define $i' = \inf\{i < s \leq t \ : \ x_s > 0\}$, and then $\sum_{s = 1}^{j-i'} \frac{x_{s+i'}}{\beta^s} = \beta^{i'-i} K \geq 1$. Furthermore, it suffices to assume $j - i \geq 2k$, since if $\sum_{s = 1}^{j-i} \frac{x_{s+i}}{\beta^s} \geq 1$ and $j - i < 2k$, then clearly 
$\sum_{s = 1}^{2k} \frac{x_{s+i}}{\beta^s} \geq \sum_{s = 1}^{j-i} \frac{x_{s+i}}{\beta^s} \geq 1$.

To this end, consider any $0 \leq i < j$ where $j - i \geq 2k$ and $x_{i+1} > 0$. In the decomposition of $y$ above into portions $w_i$ with no $0^k$ and $r_i$ runs of $0$s of length at least $k$, $x_{i+1}$ must be part of some $w_j$. Therefore, some
$s \in S$ is in the interval $(i, i+2k]$ (either $i+1$ was at the end of $w_j$, in which case $i+1 = s_j^{(1)}$ by definition, or
some $s_j^{(m)} \in (i,i+2k]$ by the fact that consecutive elements of $S$ within $w_j$ are separated by at most $2k$). By definition of 
$x$, $x(s) = 0$ and $y(s) > 0$; recall also that $x$ is coordinatewise less than $y$ on all of $\mathbb{N}$. 

This implies that 
\[
\sum_{s = 1}^{2k} \frac{x_{s+i}}{\beta^s} \leq 
\sum_{s = 1}^{2k} \frac{y_{s+i}}{\beta^s} - \beta^{-2k} \leq
\left(\frac{\beta+\delta}{\beta}\right)^{2k} \sum_{s = 1}^{2k} \frac{y_{s+i}}{(\beta+\delta)^s} - \beta^{-2k} <
\left(\frac{\beta+\delta}{\beta}\right)^{2k} - \beta^{-2k}. 
\]
(The second inequality follows from (\ref{beta2}) since $y \in X_{\beta+\delta}$.) The final expression is equal to $\frac{(\beta + \delta)^{2k} - 1}{\beta^{2k}}$, which is less than $1$ by assumption on $\delta$. We've verified (\ref{beta2}) for $x$, and so 
$x \in X_{\beta}$, completing the proof.

\end{proof}






We can now prove that the Hausdorff dimension of the limit set of $X_\beta$ is continuous as a function of $\beta$.

\begin{theorem}\label{cont}
For any CIFS $\{\phi_i\}_{i \in A}$ with $A = \{0, \ldots, j\}$, the Hausdorff dimension of $J(X_{\beta})$ is continuous
for $\beta \in [0, j+1]$.
\end{theorem}

\begin{proof}
We begin by proving a general upper bound on $\frac{Z(X_{\beta'}, n, t)}{Z(X_{\beta}, n, t)}$ for $\beta' > \beta$, which will be used to estimate $HD(J(X_{\beta'})) - HD(J(X_{\beta}))$ via the Bowen formula. For ease of notation, we define $\delta(\beta, k) = (1 + \beta^{2k})^{1/2k} - \beta$, so that the conclusion of Lemma~\ref{beta} holds when $\delta < \delta(\beta, k)$.

Choose any $k > 1$ and $\beta' < \beta + \delta(\beta,k)$. Then for every 
$y \in X_{\beta'}$, there exists $T$ with gaps greater than $k$ so that $y$ contains nonzero letters at locations in $T$, and changing those letters to $0$ yields $x \in X_{\beta}$. 

Using this, we can define, for any $n$, a map $f_n: \mathcal{L}_n(X_{\beta'}) \rightarrow \mathcal{L}_n(X_{\beta})$; to any $v \in \mathcal{L}_n(X_{\beta'})$, choose some $y \in X_{\beta'}$ beginning with $v$, then apply the replacement of Lemma~\ref{beta} to obtain $x \in X_{\beta}$, and define $f_n(v)$ to be the prefix of $x$ of length $n$. 

For any $x \in \mathcal{L}_n(X_{\beta})$, any $y \in f_n^{-1}(x)$ is completely determined by a choice of a subset $T' \subset \{1, \ldots, n\}$ with gaps greater than $k$ and the letters $y(t) \in \{1, \ldots, j\}$ for $t \in T'$. Clearly $|T'| \leq \lceil n/k \rceil$, and the number of such subsets is bounded by
\[
\sum_{i = 0}^{\lceil n/k \rceil} \binom{n}{i} \leq (\lceil n/k \rceil + 1) \binom{n}{\lceil n/k \rceil}.
\]
Since $\binom{n}{i}$ is increasing for $0 \leq i \leq n/2$, $|f_n^{-1}(x)| \leq j^{\lceil n/k \rceil} (\lceil n/k \rceil + 1) \binom{n}{\lceil n/k \rceil}$. For simplicity, we assume $n$ is a multiple of $k$, in which case $n/k$ is a positive integer. This yields
\begin{equation}\label{preimbd}
|f_n^{-1}(x)| \leq j^{n/k} (2n/k) \binom{n}{n/k}.
\end{equation}

We also require a bound relating $\left\| (\phi_w)' \right \|$ and $\left\| (\phi_{f_n(w)})' \right \|$. To see this, note that we can write
\begin{equation}\label{decomp}
w = v_1 a_1 v_2 a_2 \ldots a_{|T|} v_{|T|+1}, f_n(w) = v_1 0 v_2 0 \ldots 0 v_{|T|+1}
\end{equation}
where each $v_i$ has length at least $k$ and each $a_i \in \{1, \ldots, j\}$. For any words $u, v$ over $A$, the following inequality is a simple consequence of the bounded distortion property (iv') for shift-generated conformal iterated constructions:
\begin{equation}\label{genbd}
K^{-1} \left\| (\phi_{u})' \right\| \left\| (\phi_{v})' \right\| \leq \left\| (\phi_{uv})' \right\| \leq \left\| (\phi_{u})' \right| \left\| (\phi_{v})' \right\|
\end{equation}
(here $K$ is the constant from (iv') and is independent of the words $u,v$). Repeatedly applying this to (\ref{decomp}) yields
\begin{equation}\label{fbd}
\left\| (\phi_w)' \right \| \leq K^{|T|} \left\| (\phi_{f_n(w)})' \right \| \leq K^{n/k} \left\| (\phi_{f_n(w)})' \right \|.
\end{equation}

Combining (\ref{preimbd}) and (\ref{fbd}) yields
\begin{multline}
Z(X_{\beta'},n,t) = \sum_{w \in \mathcal{L}_n(X_{\beta'})} \left\| (\phi_w)' \right \|^t
\leq K^{2nt/k} \sum_{w \in \mathcal{L}_n(X_{\beta'})} \left\| (\phi_{f_n(w)})' \right \|^t\\
\leq (jK^{t})^{n/k} (2n/k) \binom{n}{n/k} \sum_{v \in \mathcal{L}_n(X_{\beta})} \left\| (\phi_{v})' \right\|^t
= (jK^{t})^{n/k} (2n/k) \binom{n}{n/k} Z(X_{\beta},n,t).
\end{multline}

We now take logarithms, divide by $n$, and let $n$ approach infinity (technically along the subset of multiples of $k$) to get
\begin{multline}\label{mainbd}
P(X_{\beta'}, t) \leq P(X_\beta, t) + \frac{\ln j - t \ln (K^{-1})}{k} - k^{-1} \ln(k^{-1}) - (1 - k^{-1}) \ln (1 - k^{-1})\\
\leq P(X_\beta, t) + \frac{\ln j - t \ln (K^{-1}) + 2 \ln k}{k}
\end{multline}
(recalling that $k > 1$, and so $- (1 - k^{-1}) \ln (1 - k^{-1}) \leq - k^{-1} \ln(k^{-1}) = \frac{\ln k}{k}$).

Finally, we fix any $\beta \in [0, j+1]$ and define $h = HD(J(X_{\beta}))$; by Proposition~\ref{Bowenf}, $P(X_{\beta}, h) = 0$. Take any sequence $\beta_n$ approaching $\beta$ from above and define $h_n = HD(J(X_{\beta_n}))$; again by Proposition~\ref{Bowenf},
$P(X_{\beta_n}, h_n) = 0$. For any $\epsilon > 0$, choose $k$ so that
$\frac{\ln j - HD(M) \ln (K^{-1}) + 2 \ln k}{k} < \epsilon$. Then for large enough $n$, $\beta_n < \beta + \delta(\beta, k)$, and so, using (\ref{mainbd}) (and noting that all $h_n$ are bounded from above by $HD(M)$),
\[
0 = P(X_{\beta_n}, h_n) \leq P(X_{\beta}, h_n) + \frac{\ln j - h_n \ln (K^{-1}) + 2 \ln k}{k} < P(X_{\beta}, h_n) + \epsilon \leq P(X_{\beta}, h) + \epsilon = \epsilon.
\]
This means that $-\epsilon < P(X_{\beta}, h_n) \leq 0$ for sufficiently large $n$. Since $\epsilon$ was arbitrary and
$P(X_{\beta}, t)$ is a monotone function of $t$ with a unique root, this implies that $h_n \rightarrow h$ and so that
$HD(J(X_{\beta}))$ is right-continuous.

For left-continuity, assume that $\beta_n$ approaches $\beta$ from below, choose arbitrary $\epsilon > 0$, and define $k, h, h_n$ as before. There is a tiny subtlety; we need $\beta < \beta_n + \delta(\beta_n, k)$, which theoretically could be an issue since
$\delta(\beta_n, k)$ is not constant. However, we note that $\delta(x,k)$ is monotone decreasing in $x$, and so it suffices to satisfy
$\beta < \beta_n + \delta(\beta, k)$, which is obviously satisfied for sufficiently large $n$. We then again use (\ref{mainbd}):
\[
0 = P(\beta, h) \leq P(X_{\beta}, h_n) \leq P(X_{\beta_n}, h_n) + \frac{\ln j - h_n \ln (K^{-1}) + 2 \ln k}{k} < P(X_{\beta_n}, h_n) + \epsilon = \epsilon.
\]
This means that $0 \leq P(X_{\beta}, h_n) < \epsilon$ for sufficiently large $n$. As above, this implies $h_n \rightarrow h$ and so that $HD(J(X_{\beta}))$ is left-continuous, completing the proof.

\end{proof}

Theorem~\ref{cont}, along with the Bowen formula, allows us to prove that many classes of CIFS and CDGMS have maximal extended Hausdorff dimension spectrum.

\begin{theorem}\label{spectrum}
For any CIFS $\{\phi_i\}_{i \in A}$ on $Y$ over a finite alphabet $A$, the extended Hausdorff dimension spectrum
$\{HD(J(X)) \ : \ X \subset A^{\mathbb{N}}\}$ is $[0, HD(J(Y))]$.
\end{theorem}

\begin{proof}
Without loss of generality, we assume $A = \{0, \ldots, j\}$.
The shift $X_0 = \{0^{\infty}\}$ corresponding to $\beta = 0$ is a singleton, so $HD(J(X_0)) = 0$. Similarly, the shift $X_{n+1}$ corresponding to $\beta = j + 1$ is the full shift $\{0, \ldots, j\}^{\mathbb{N}}$, so $J(X_{j + 1}) = J(Y)$ and
$HD(J(X_{j+1})) = HD(J(Y))$. Theorem~\ref{spectrum} now follows from Theorem~\ref{cont} by the Intermediate Value Theorem.
\end{proof}




\begin{theorem}\label{spectrum2}
For any CGDMS $\{\phi_i\}_{i \in A}$ over a finite alphabet $A$ defined by a topological Markov chain 
$Z \subset A^\mathbb{N}$, the restricted extended Hausdorff dimension spectrum $\{HD(J(X)) \ : \ X \subset Z\}$ is $[0, HD(J(Z))]$.
\end{theorem}

\begin{proof}
We first look at the scenario in which $Z$ is irreducible. In this case for all $i, j \in A$
there is a word $w_{i,j}$ for which $i w_{i,j} j \in \mathcal{L}(Z)$. 


Fix any letter $a \in A$. For every $m$, we will define a subshift $Z_{m} \subset Z$ which behaves like a full shift with finite alphabet. Specifically, for each  $v \in \mathcal{L}_m(Z)$, there exist $u_v, u'_v$ of length less than $N$ for which $a u_v v u'_v a \in \mathcal{L}(Z)$ (just define $u_v = w_{a, v_1}$ and $u'_v = w_{v_{|v|}, a}$). For any $v \in \mathcal{L}_m(Y)$, 
define $t_v := a u_v v u'_v$, and define $Z_{m}$ to be the set of all shifts of infinite concatenations of words in the set $\{t_v\}$. Since $Z$ is a topological Markov chain, $Z_{m} \subset Z$. Fix an arbitrary enumeration $\{v_i\}_{i=1}^M$ of $\mathcal{L}_m(Z)$, and for any $\beta$-shift $X_{\beta}$ with $1 < \beta \leq M + 1$, define $Z_{m,\beta}$ to be the subshift consisting of all shifts of concatenations of the form 
$t_{v_{n_1}} t_{v_{n_2}} t_{v_{n_3}} \ldots$ where $n_1 n_2 n_3 \ldots \in X_{\beta}$. 

It's not difficult to check that the collection $\{\phi_{t_{v_i}}\}_{1 \leq i \leq M}$ is a conformal IFS . For any $\beta$,
define $J_{\beta} = J(Z_{m, \beta})$ the limit set for the shift-generated conformal construction from the original IFS and
$J'_{\beta} = J(X_{\beta})$ the limit set for the shift-generated conformal construction from the IFS $\{\phi_{t_{v_{n_i}}}\}$. These sets need not be equal, but they are related; $J'_{\beta}$ corresponds to intersections of images of compositions for only those right-infinite sequences 
in $Z_{m, \beta}$ which are infinite concatenations of the form $t_{v_{n_1}} t_{v_{n_2}} \ldots$ (as opposed to a shift of such a sequence). Therefore, clearly $J'_{\beta} \subset J_{\beta}$, so $HD(J'_{\beta}) \leq HD(J_{\beta})$. In addition, since very point in $Z_{m,\beta}$ can be written as 
$s t_{v_{n_1}} t_{v_{n_2}} \ldots$ for some suffix $s$ of a word $t_v$, $J_\beta$ is a subset of the union of all images of $J'_\beta$ under $\phi_s$ for such suffixes $s$. Then, since all such $\phi_s$ are Lipschitz and Hausdorff dimension is nonincreasing under Lipschitz maps, each such image has Hausdorff dimension less than or equal to $HD(J'_\beta)$, meaning that $HD(J_\beta) \leq HD(J'_\beta)$. So, $J_\beta$ and $J'_\beta$ have equal Hausdorff dimension. By Theorem~\ref{cont}, $HD(J'_{\beta})$ is continuous as a function of $\beta$, and so the same is true of $HD(J_\beta) = HD(J(Z_{m, \beta}))$. 


Clearly $Z_{m,0}$ consists of a single periodic orbit (corresponding to $t_{v_0})$ (implying $HD(J(Z_{m,0})) = 0$) and $Z_{m,M+1} = Z_m$. So for every $m$, by the Intermediate Value Theorem, limit sets of shift-generated conformal constructions for subshifts of $Z_m$ achieve all Hausdorff dimensions in $[0, HD(J(Z_m))]$.

We now estimate $HD(J(Z_{m,M+1})) = HD(J(Z_{m}))$ from below using Proposition~\ref{Bowenf}. As a preliminary, we note that by (\ref{genbd}), for any $v_1, \ldots, v_k \in \mathcal{L}_m(Z)$, 
\begin{equation}\label{bd1}
\left\| (\phi_{t_{v_1} \ldots t_{v_k}})' \right\| \geq K^{-3k} \prod_{i = 1}^k \left\| (\phi_{a})' \right\| \left\| (\phi_{u_{v_i}})' \right\| \left\| (\phi_{v_i})' \right\| \left\| (\phi_{u'_{v_i}})' \right\| \geq K^{-3k} L^{2k} W^k \prod_{i = 1}^k \left\| (\phi_{v_i})' \right\|,
\end{equation}
where $W = \left\|(\phi_a)'\right\|$ and $L = \min_{i, j \in A} \left\|(\phi_{w_{i,j}})'\right\|$ (recalling that all $u_{v_i}$ and $u'_{v_i}$ are some $w_{i,j}$).

Define $N = \max(|w_{i,j}|)$; then all $u_{v_i}$ and $u'_{v_i}$ have length between $1$ and $N$, and for any $k > 0$, all concatenations
$t_{v_1} \ldots t_{v_k}$ have length between $k(m+1)$ and $k(m+2N+1)$. Therefore, by (\ref{bd1}), for any $t$,  
\begin{multline}\label{bd2}
\sum_{n = k(m+1)}^{k(m+2N+1)} Z(Z_{m}, n, t) = 
\sum_{n = k(m+1)}^{k(m+2N+1)} \sum_{s \in \mathcal{L}_{n}(Z_{m})} \left\| (\phi_w)' \right\|^t \geq
\sum_{v_{1}, \ldots, v_{k} \in \mathcal{L}_m(Z)} \left\| (\phi_{t_{v_1} \ldots t_{v_k}})' \right\|^t \\
\geq (K^{-3} L^2 W)^{kt} \sum_{v_{1}, \ldots, v_{k} \in \mathcal{L}_m(Z)} \prod_{i = 1}^k \left\| (\phi_{v_i})' \right\|^t
= (K^{-3} L^2 W)^{kt} (Z(Z, m, t))^k.
\end{multline}
Let's specifically define $t_m = HD(J(Z_{m}))$; then by Proposition~\ref{Bowenf}, $P(Z_m, t_m) = 0$. Therefore, the quantities 
$Z(Z_{m}, n, t_m)$ are subexponential in $n$. This means that when we take logarithms, divide by $km$, and let $m \rightarrow \infty$ in 
(\ref{bd2}), we get
\[
0 \geq \frac{h_m \ln(K^{-3} L^2 W)}{m} + P(Z, h_m) \Longrightarrow P(Z, h_m) \leq \frac{-\ln(K^{-3} L^2 W) h_m}{m} \leq \frac{C}{m},
\]
where $C = -\ln(K^{-3} L^2 W) HD(M)$. In addition, since $h_m \leq h = HD(J(Z))$, $P(Z, h) = 0$ by Proposition~\ref{Bowenf}, and $P(Z, t)$ is decreasing, $P(Z, h_m) \geq 0$. Therefore, for all $m$, $P(Z, h_m) \in [0, \frac{C}{m}]$, meaning that $P(Z, h_m) \rightarrow 0$. Since $P(Z, t)$ is a monotone function of $t$ with a unique root $h$, $HD(J(Z_m)) = h_m \rightarrow h = HD(J(Z))$.
Since we already showed that $\{HD(J(X)) \ : \ X \subset Z\}$ contains $[0, HD(J(Z_m))]$ for all $m$, it also contains $[0, HD(J(Z)))$. Since clearly $HD(J(Z))$ is achieved by $X = Z$, this completes the proof for the case when $Z$ is irreducible.

Let us now move to the more general case, when $Z$ is not assumed to be irreducible. In this instance, the finite alphabet $A$ has a unique decomposition $A = A_1 \cup A_2 \cup \ldots A_l$, where each restriction $Z \cap (A_k)^\mathbb{Z}$ is an irreducible topological Markov chain over $A_k$ and each $A_k$ is inclusion maximal with respect to this restriction. Of these sets, one of them, say $A_j$, generates a limit set whose Hausdorff dimension is $HD(J(Z))$ (see \cite{GMR} for more details; in particular Theorem 3.5 and Corollary 3.6). We will denote by $Z_j$ the intersection $Z \cap (A_j)^\mathbb{Z}$.

Given that $Z_j$ is irreducible,  the restricted extended Hausdorff dimension spectrum $\{HD(J(X)) \ : \ X \subset Z_j\}$ is $[0, HD(J(Z_j))]$. Since $Z_j \subset Z$ and $HD(J(Z_j)) = HD(J(Z))$, this finishes the proof.

\end{proof}





\begin{theorem}\label{spectrum3}
For any CGDMS $\{\phi_i\}_{i \in A}$ over a countably infinite alphabet $A$ defined by a finitely irreducible topological Markov chain $Z \subset A^\mathbb{N}$, the restricted extended Hausdorff dimension spectrum $\{HD(J(X)) \ : \ X \subset Z\}$ is $[0, HD(J(Z))]$. 
\end{theorem}

\begin{proof}
Without loss of generality, let's assume $A = \{0, 1, \ldots\}$. Then for any finite subset $B \subset A$ we know from the above
that the extended Hausdorff dimension spectrum contains all reals in $[0, HD(J(B^\mathbb{N}))]$. But it is known (see Theorem 3.15 in \cite{MU}) that
$HD(J(Z)) = \sup_B HD(J(B^{\mathbb{N}} \cap Z))$, and so the extended Hausdorff dimension spectrum contains $[0, HD(J(Z)))$. 
The proof is complete.
\end{proof}

\end{document}